\documentclass[preprint]{article}
\usepackage{hyperref}
\usepackage{amsmath}
\usepackage{amsfonts}
\usepackage{amsthm}
\usepackage{amssymb}
\usepackage{graphicx}
\usepackage{float}
\usepackage{caption}
\usepackage{subcaption}
\usepackage{wrapfig}

\usepackage{amssymb,amsmath,graphicx,psfrag,amsfonts,latexsym,amsthm,amscd}
\usepackage[all, knot]{xy}
\xyoption{arc}

\newtheorem{theorem}{Theorem}[section]
\newtheorem{lemma}[theorem]{Lemma}
\newtheorem{corollary}[theorem]{Corollary}

\newtheorem{prop}[theorem]{Proposition}
\theoremstyle{definition}

\begin{document}
\title{A spatial version of Tutte's conflict graph}

\author{Joel Foisy}

\maketitle

\begin{abstract}
Tutte showed that a graph $G$ is planar if and only if the conflict graph associated to every cycle of $G$ is bipartite. Inspired by this, we define a (not necessarily unique) signed conflict graph associated to a maximally planar subgraph of a nonplanar graph such that if a nonplanar graph $G$ has a flat embedding, there is a balanced conflict graph associated to every maximally planar subgraph of $G$. We conjecture that a graph $G$ is intrinsically linked if and only if every maximal planar subgraph of $G$ has every possible conflict graph unbalanced. Foisy and Raimondi's finding of every maximal planar subgraph of every graph in the Petersen family having every possible conflict graph unbalanced provides further supporting evidence for this conjecture.
\end{abstract}

\renewcommand{\thefootnote}{\fnsymbol{footnote}} 
\renewcommand{\thefootnote}{\arabic{footnote}} 
\section{Introduction}

Tutte showed that a graph $G$ is planar if and only if every cycle of $G$ has a bipartite conflict graph \cite{tutte}. This can be proven using Kuratowki's theorem (see, for example, \cite{bl}), and it is also the key idea in Dirac and Schuster's proof of Kuratowski's Theorem \cite{DS2}, \cite{DS}. For a graph $G$, the {\em conflict graph} (also called {\em overlap graph} by some authors, see for example \cite {bl}) of a cycle $C$ in $G$ is the graph in which each vertex represents a fragment of $C$ and two vertices are adjacent if their corresponding fragments conflict. A {\em fragment} of $C$ (also called a {\em bridge} of $C$ by some authors, see again for example \cite{bl}) are either connected components of $G-C$, along with their edges of attachment to $C$, or chords of $C$. Two $C$-fragments, $A,B$ \textit{conflict} if they have three common vertices of attachment to $C$ or if there are four vertices, $v_1, v_2, v_3, v_4$, in cyclic order on $C$, such that $v_1$ and $v_3$ are vertices of attachment for $A$, and $v_2$ and $v_4$ are vertices of attachment for $B$. In using the terms fragment and conflict graph, we are following the conventions used in West's book \cite{west}.

\begin{figure}
\hskip 1.15in \includegraphics[scale=.185]{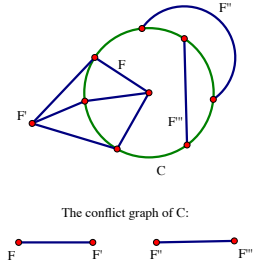}
\caption{A planar graph, and the bipartite conflict graph associated to a cycle $C$.}
\end{figure}
\vskip .01in

In Tutte's and Dirac-Schuster's work on the conflict graph, given a bipartite conflict graph, one can start with an embedded cycle $C$, and build a planar embedding by placing fragments either inside or outside of the cycle $C$, placing fragments that correspond to adjacent vertices in the conflict graph (conflicting fragments) on opposite sides of the cycle. Here we propose a generalization of the conflict graph for spatial graphs, where the role of cycle will be played by maximally planar subgraphs. Given a connected nonplanar graph $G$, one starts with an embedding of a maximally planar subgraph $M$ on $S^2$, sitting nicely (say, as a unit sphere) in 3-space ($M$ \textit{maximally planar} here means $M$ is planar, but including any edges of $G-M$ results in a nonplanar graph), and then places $M$-fragments (edges in $G-M$) either inside or outside the sphere, with conflicting fragments on different sides, and anti-conflicting fragments on the same side. The hope is that when the conflict graph is balanced (a generalization of bipartite for signed graphs), this placement of fragments (with a possible adjustment of the embedding among fragments on the same side) will lead to a flat embedding. (A \textit {flat} spatial embedding of $G$ is an embedding in which it is possible to attach a disk to every cycle of $G$ such that the interior of the disk is disjoint from the graph--$G$ having a flat spatial embedding is equivalent to $G$ having a linkless spatial embedding \cite{robertson seymour thomas}).

We give a definition of the (signed) conflict graph in this situation and show that if $G$ has a \textit{spherically flat embedding} (that is, a flat spatial embedding for which every planar subgraphs lies on a sphere that intersects $G$ only in the planar subgraph), then an associated conflict graph is balanced for every maximally planar subgraph of $G$. We further show in Section 2 that every flat embedding is in fact a spherically flat embedding. This generalizes a result of Wu \cite{wu} that every flat spatial embedding of a planar graph is \textit {spherical} (there exists a topological 2-sphere that contains the flatly embedded planar graph).



We further conjecture that $G$ has a flat embedding precisely when every possible conflict graph is balanced for every maximally planar subgraph $M$ (and further, if one maximally planar subgraph of $G$ has a balanced (unbalanced) conflict graph, then all do). Similarly, we conjecture that having every possible conflict graph balanced for every maximally planar subgraph $M$ is closed under minors. This, together with Foisy and Raimondi's finding \cite{fr} of every possible conflict graph unbalanced for every maximally planar subgraph of every graph in the Petersen family, and Robertson, Seymour and Thomas' work \cite{robertson seymour thomas}, would establish that every maximally planar subgraph of an intrinsically linked graph has every possible conflict graph unbalanced.\\

In terms of the layout of the paper after Section 2, in Section 3, we define what it means, given a nonplanar graph $G$ and maximal planar subgraph $M$ embedded in $S^2$, for two $M$-fragments $F$ and $F'$ to strongly anticonflict. We show in Proposition \ref{linked} that if $F$ and $F'$ strongly anticonflict, then by placing $F$ and $F'$ on the same side of the sphere containing $M$ (which we are considering nicely sitting inside 3-space), then the resulting embedding always contains a pair of nonsplit linked cycles. 

In Theorem \ref{noconflict}, we show that if $F$ and $F'$ do not strongly conflict, then we can extend the embedding of $M$ to a flat embedding of $M \cup F \cup F'$ that has $F$ and $F'$ on the same side of $S^2$. We later define strongly anticonflicting fragments and prove analogous results for them. For example, in Proposition \ref{noconflict2}, we show that for not strongly anticonflicting fragments $F$ and $F'$, then a flat embedding of $M \cup F \cup F'$ exists, with $F$ and $F'$ placed on opposite sides of the $S^2$ that contains $M$.

In the last section, starting with an embedding of a maximally planar subgraph $M$ in $S^2$, we use the placement of strongly (anti)conflicting fragments on appropriate sides of the sphere containing $M$--extending the planar embedding of $M$ in $S^2$ to a possibly flat embedding of $G$ into space, to help define implicitly conflicting fragments. From there, we define the (signed) conflict graph of an $S^2$ embedding of M as the signed graph having the fragments of M as vertices and including negative edges between strongly or implicitly conflicting fragments, and positive edges between strongly or implicitly anti-conflicting fragments. \\
 

 We end this introductory section with a few more definitions. We say a graph is {\em intrinsically linked} if, in every spatial embedding, there is a pair of cycles that form a nonsplit link. A graph $H$ is a \emph{minor} of a graph $G$ if $H$ can be obtained from $G$ by a sequence of vertex deletions, edge deletions and edge contractions. Recall that the seven Petersen Family Graphs are the complete minor-minimal set of intrinsically linked graphs \cite{robertson seymour thomas}, and consist of $K_6$ and the six graphs derivable from it by $\Delta-Y$ and $Y-\Delta$ exchanges, denoted by $K_{3,3,1}, P_7, P_8, P_9, K_{4,4}-e$ and $P_{10}$, where $P_{10}$ is the classic Petersen Graph.    


We work in the PL category.

\section{Flat implies spherically flat}
In this section, we prove some results that will allow us to show that, for a flat embedding of a graph $G$, every planar subgraph $P$ of $G$ lies in an embedded sphere that meets $G$ only in $P$. We will make use of B\"ohme's Lemma:
\begin{lemma} (B\"ohme) \cite{bohme}: Let $\phi$ be a flat embedding of $G$ and let $C_1, ..., C_n$ be a family of cycles of $G$ such that for every $i \neq j$, the intersection of $C_i$ and $C_j$ is either connected or null. Then there exist pairwise disjoint open disks $D_1, ..., D_n$, disjoint from $\phi(G)$ and such that $\phi(C_i)$ is the boundary of $D_i$ for $i=1, 2, ..., n$.
\label{bohme}
\end{lemma}


Now we introduce another Lemma, which roughly says that for a graph $G$ embedded flatly, and $C$ a cycle in $G$, we can find multiple pair-wise inequivalent paneling disks for $C$. If the cycle is separating, we can find at least 3 such paneling disks. For the following result, for a positive integer $n$, take $[n]=\{1,2,...,n\}$.

\begin{lemma}  Let $G$ be a connected graph with a flat embedding, $\phi(G)$, and $C$ a cycle of $G$.  Let $\{F_i, i \in [n]\}$ be the collection of all $C$-fragments in $G$.  Then for every flat embedding of $G$, there is an ordering of the $\{F_i\}$ and a collection of $3-$balls, $\{B_i : i \in [n]\}$ such that 
\begin{enumerate}
\item  For every $i \in [n]$, $\partial(B_i)=b_{i_1} \cup b_{i_2}$, where each $b_{i_j}$ is a $2-$cell with $\partial(b_{i_j})=\phi(C)$.  That is, each $b_{i_j}$ is a paneling disk for $\phi(C)$.
\item  For every $i \in [n]$, $F_i \subset B_i$
\item  For $i, j \in [n]$, with $j > i+1$, $B_i \cap B_j = \phi(C)$.

\item For $i, j=i+1 \in [n]$, $B_i \cap B_{j} =b_{i_2}=b_{j_1}$.
\item  For every $i \in [n]$, $B_i \cap \phi(G) = \phi(F_i) \cup \phi(C)$.
\end{enumerate}
\label{magic}
\end{lemma}

\begin{proof}  
Assume the hypotheses. We will prove the result by induction on the number of vertices of $G$, where the number of vertices is at least $4$ (otherwise, $G$ is acyclic or just a 3-cycle).  For our base case, we note that, by a result of Wu \cite{wu}, every flat embedding of a planar graph is spherical (that is, there is a topological $2-$sphere in ${\mathbb R}^3$, such that the embedded graph lies on the $2-$sphere). The result can be seen for every graph on 4 vertices. A nonseparating cycle, with a paneling disk, can be extended to a ball that contains the graph minus the cycle in its interior, with the paneling cycle as part of the boundary. The ball contains another paneling disk for the cycle, as described in the theorem. A separating cycle in a spherical graph similarly satisfies the result, by merely deforming the fragments slightly, with a paneling disk in the sphere, another panel pushed off the first in direction of one fragment, and the third panel pushed off the first in the direction of the other fragment.

Now suppose the claim is true for graphs with $n$ vertices.  Let $G$ be a connected graph with $n+1$ vertices that has a flat embedding.  Flatly embed $G$. Let $C$ denote a cycle of $G$. If possible, we select our edge to contract, e, to be contained in a $C$-fragment. By induction, we can find $B_i$ meeting the conclusions of the lemma. Then the vertex splitting used to get $e$ back can be done preserving the fragment and all of the desired properties. Otherwise, there is no edge entirely contained in a fragment, that is to say, every edge that is not contained in $C$ has an endpoint in $C$.  Thus $G$ is a subgraph of $C + \bar{K}_n$, which is a wheel with $n$ hubs, where some spokes may be missing.
If $C$ has more than 3 vertices, we can contract an edge, $e$, of the  cycle.  The embedding of $G$ induces an embedding of $H=G/e$, with corresponding cycle $C/e$. Note that contracting edge $e$ will neither create nor destroy $C-$fragments. By the induction hypothesis, every fragment of $C/e$ in $H$, $H_i$, lies in a $3-$ball, $B_i$, meeting all of the assumptions spelled out above.  We then expand back the vertex corresponding to $e$ to get the edge $e$, and to get the original embedding of $G$ back with each $C-$fragment meeting the conclusion of the lemma.  

Otherwise, the cycle $C$ has exactly $3$ vertices and there is no edge entirely contained in a fragment (an endpoint of every edge must lie in $C$). Thus $G$ is a subgraph of $K_{1,1,1,n}$. For a flat embedding of $G$, we extend to a flat embedding of $K_{1,1,1,n}$. The lemma applies in this case, as the collection of all 3-cycles satisfies the hypotheses of B\"ohme's Lemma. These cycles can be simultaneously paneled, by Lemma \ref{bohme}, and this paneling forms the $B_i$ as described. This concludes our proof.

\end{proof}

\begin{figure}[H]

\hskip 1in  \includegraphics[scale=.21]{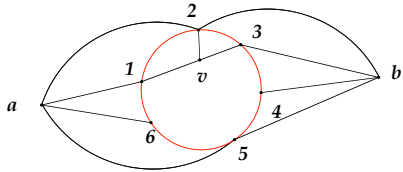}

\vskip -.35in
\caption{Here the cycle $(1,2,3,4,5,6)$ has $3$ fragments, and the reader can check that there are at least four different in-equivalent paneling disks for the cycle.}
\label{example}
\end{figure}

We observe that if $G$ has a flat embedding $\phi(G)$, and $C$ is a cycle in $G$, then the previous Lemma establishes that it is possible to choose equivalence class representatives for some classes of paneling disks for $C$ such that each $C-$fragment is sandwiched between exactly two such disks (meaning the two disks form a sphere $S$ that contains the fragment in the bounded component of ${\mathbb R}^3-S$), and there is an injection between a collection of such disk representatives and all $C$ fragments plus 1.  There can be more disk classes than $C$-fragments, see Figure \ref{example}. 

\vskip -.4in

The next result is a key step in showing that that every maximal planar subgraph of every nonplanar flat graph has a balanced conflict graph.

\begin{theorem} If $G$ is nonplanar and has a flat embedding, and if $P$ is a 2-connected spanning planar subgraph of $G$, then for a given flat embedding of $G$, there exists a sphere $S$ that contains the embedding of $P$, and $S$ is disjoint from the interiors of all edges not in $P$ (that is, all $P-$fragments).
\label{disjoint}
\end{theorem}

\noindent
\begin{proof}  First, as a warm-up, let $P$ be a 3-connected spanning planar subgraph of $G$, where $G$ is embedded flatly. As $P$ is 3-connected, then as Robertson, Seymour and Thomas point out in \cite{robertson seymour thomas} (Proof of (2.2)), the face boundaries of $P$ meet the hypotheses of Lemma \ref{bohme}, then $\phi(P) \cup D_1 \cup... \cup D_n$ form the desired sphere.\\

Now suppose $P$ is a 2-connected spanning planar subgraph of $G$, where $G$ is embedded flatly. Let $S$ be a sphere that contains $P$, as guaranteed by Wu's result. Choose $S$ so that it intersects $G-P$ transversely, and so that it has a minimal number of intersections with $G-P$. If there are no intersections, we are done. 
Otherwise, let $C$ be the boundary of a face of $P$ in $S$ that intersects $G-P$. We note that the boundary is a cycle because $P$ is 2-connected, see, for example \cite{west}. In particular, let $e$ be an edge of $G-P$ that intersects $S$, orient $e$, and take $C$ to be the cycle bounding the last region that $e$ intersects (see Figure \ref{fig:pierced}). Call the endpoint of $e$ following this intersection point, $v$ (see Figure \ref{fig:lastint}). We may assume that that $v$ does not lie in the region bounded by $C$, nor on $C$, for then we could reduce the number of intersections of $e$ and $S$ by a simple isotopy of $e$.

\vskip -.3in
\begin{figure}[H]
\vskip .1in
\hskip 1.45in
\includegraphics[scale=.14]{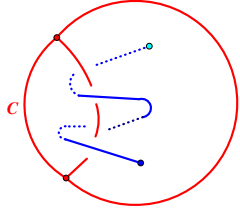}
\caption{The cycle $C$ bounds the last region pierced by $e$.}
\label{fig:pierced}
\end{figure}

Take a paneling disk, $D$, for $C$ that intersects $S$ transversally and with minimal components. The intersection of $D$ and $S$ will be the union of $C$ and circles and arcs with endpoints in $C$ (see Figure \ref{fig:intersections}).
\vskip -.1in
\begin{figure}[H]
\vskip .1in
\hskip 1.45in
 \includegraphics[scale=.14]{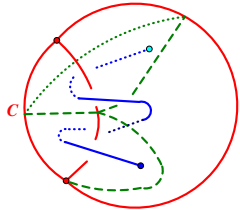}
 \caption{The green curves represent the intersection of the paneling disk for $C$ with $S$.}
 \label{fig:intersections}
\end{figure}

The endpoint $v$ will be contained in a component of $S-(C \cup D)$, which we will denote $N_v$. Thanks to Lemma \ref{magic}, we may assume the paneling disk for $C$, together with the face determined by $C$, bounds a topological ball that contains $v$ and a segment of $e$. The disk comes in underneath $N_v$. How does Lemma \ref{magic} apply? It is possible that the paneling disk, $D'$ came into $C$ overlapping the face of $C$ near the vertices of attachment of $N_v$.  By Lemma \ref{magic}, we can find another a paneling disk $D$ such that the $C-$fragment containing $v$ is contained in the ball bounded by $D \cup D'$. 

\vskip .05in
\begin{figure}
\vskip .18in
\hskip1.45in
\includegraphics[scale=.14]{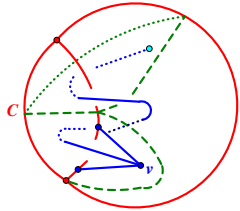}
\caption{The endpoint $v$ lies in a component of $S-(C \cup D)$. The component may be more complicated than shown here.}
\label{fig:lastint}
\end{figure}

\begin{figure}
 \hskip 1.45in
\includegraphics[scale=.14]{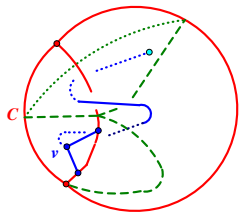}
\caption{After the ambient isotopy, the number of intersections of $G-P$ with $S$ has decreased. If more than one fragment terminated in $N_v$, the number of intersections may have decreased by more than one. }
\end{figure}

Let $f$ be an arbitrary $P-$fragment (edge not in $P$) between $N_v$ and another part of $P$. Orient $f$ so that it begins with the endpoint in $N_v$. Then $f$ will intersect $S$ first in the region bounded by $C$. Note there may be more than one such fragment, but the argument will still hold: all fragments will intersect $S$ first in the region bounded by $C$.\\

It follows that the portion of $P$ that lies in $N_v$ can be ambient isotoped so that it lies in the face bounded by $C$ (see Figures \ref{fig:intersections} and \ref{fig:lastint}). Note that, in terms of the embedding of $P$ on $S$, the result of the ambient isotopy is to rotate the region $R$ inside the inner most arc with endpoints on $C$, such that $C$ with this arc contains $N_v$ so that $R$ lies on the other component of $S -C$.

The result is an embedding of $P$ on $S$ that has fewer intersections with $G-P$. This is a contradiction of the minimality of $S$. This completes our proof for when $P$ is 2-connected.

\end{proof}

We prove the following that generalizes Theorem \ref{disjoint}.
\begin{theorem} If $G$ is nonplanar with a flat embedding, and if $P$ is a spanning planar subgraph of $G$, then for a given flat embedding of $G$, there exists a sphere $S$ that contains the embedding of $P$, and $S$ is disjoint from the interiors of all edges not in $P$.
\label{flat}
\label{conjecture}
\end{theorem}

\begin{proof}
First we prove the result with $P$ is 1-connected. We do so by induction on the number of blocks of $P$. Recall that a \emph{block} of a graph is a maximal connected subgraph without a cut vertex. A block is either an edge or a 2-connected subgraph. For a flat embedding of $G$, consider two blocks of $P$ connected at $v$, call them $B_1$ and $B_2$. If either $B_i$ is an edge, we can put that $B_i$ in a sphere, $S_i$, that intersects $G$ only in $B_i$. Otherwise, by Theorem \ref{disjoint}, again $B_i$ is contained in a sphere, $S_i$, that does not intersect the rest of $G$. Make the intersections between the $S_i$ transverse: thus the intersections will be circles (see, for example: \cite{bohme} or \cite{sphere}). We can so in such a way that $S_i$ intersects $G$ only on $B_i$. By the Jordan-Brouwer Theorem, $S_i$, $i=1,2$, cuts $3$-space into two components. By Alexander's result in the PL category, the bounded side of $S_i$ is a ball \cite{alex}. 

Fill in one of the spheres, say $S_2$, to obtain a ball $D_2$ with the boundary of $D_2$ to be $S_2$. By the Jordan-Brouwer Theorem, $S_1$ cuts $D_2$ into components: some on one side of $S_1$, and the others on the other side of $S_1$. We can form a tree (see, for example, \cite{sphere}) by placing a vertex in each such component, and connecting components that share a common boundary in $S_1$. Because $D_2$ is contractible to a point, the graph must be a tree. Moreover, each component must be a ball, by the argument in the next paragraph. Finally, $B_2$, being connected, can only lie in one such component. We may remove the others to obtain $D_2'$ with boundary $S_2'$ and $B_2$ in $S_2'$. Any intersections between $S_1$ and $S_2'$ can be removed, except $v$. From here, we make a third sphere (by inserting a small collar around $v$ that connects $S_1$ and $S_2'$ with small disks, disjoint from $P$, removed near $v$) that contains $B_1$ and $B_2$, and the sphere is disjoint from the rest of the graph. The general induction is nearly identical to the base step.\\

Why is each piece a ball? We show this by induction on the number of components--for transverse intersection. This is clear for two components. Suppose true for n components. Consider an intersection with n+1 components. Form the dual tree. Go to a leaf. Its corresponding component must intersect $S_1$ in just one component. Remove it and slide the intersection of the ball with $S_1$ off. By induction, the result is a ball cut into $n$ pieces, each a ball. We can then place the leaf component back, and we have a ball cut into $n+1$ balls.\\

Finally, if $P$ is disconnected, we do an induction argument on the number of different components. For the base case if we have two disjoint planar subgraphs embedded on spheres, and if the spheres intersect, we can apply a similar process as was applied above to make the spheres disjoint. Then, we can join them together by a small tube that is disjoint from the graph, and remove the interiors of the circles of attachment to make a new sphere that contains both components. The general induction is similar to the base case. This concludes our proof.
\end{proof}

\section{Strongly (anti-)conflicting fragments}

Let $G$ be a nonplanar graph with a maximally planar subgraph $M$, embedded into $S^2$. Two $M-$ fragments (edges), $F$ and $F'$, \textit{strongly conflict} if the vertices of attachment of $F$, $v_1$ and $v_2$, and the vertices of attachment of $F'$ form the 4-partition of a $K_{4,2}$ subgraph [or expansion of $K_{4,2}$] of $M$, such that, $v_1$ and $v_2$ do not lie in the same face of the induced embedding of the $K_{4,2}$. 

\begin{figure}
\vskip -.4in
\hskip 1.5in 
\includegraphics[scale=.18]{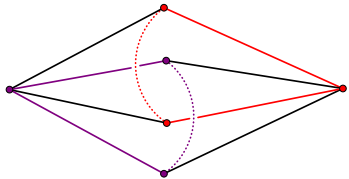}
\vskip -.23in
\caption{The dotted edges represent $M$-fragments that strongly conflict.}
\label{fig:conflict}
\end{figure}

See Figure \ref{fig:conflict} for an illustration of strongly conflicting $M$-fragments, where some edges of the subgraph $M$ may not be shown.

\begin{prop}

Let $G$ be a nonplanar graph embedded in space with (maximally) planar subgraph $M$ embedded on a sphere that meets $G$ only in $M$, and strongly conflicting $M$-fragments $F$ and $F'$. If $F$ and $F'$ lie on the same side of the sphere containing $M$, then this embedding of $G$ contains a pair of linked cycles.
\label{linked}
\end{prop}
\begin{proof}
The proof is nearly the same as Sachs' \cite{sachs} and Conway-Gordon's \cite{conwaygordon} proof that $K_6$ is intrinsically linked. Consider one such embedding. Consider disjoint cycle pairs that contain $F$ and $F'$ and edges from the $K_{4,2}$ subgraph of $M$ that force $F$ and $F'$ to strongly conflict. The sum of the linking numbers on all such pairs of cycles is one for the embedding in Figure \ref{fig:conflict}. To obtain another such embedding, one can use use ambient isotopy and crossing changes only between $F$ and $F'$. Such changes will not change the parity of the sum of the linking numbers. Thus, there exists a non-split link in every such embedding.
\end{proof}

We prove a helpful converse to the previous result. 

\begin{theorem}

Let $G$ be a connected nonplanar graph with maximally planar subgraph $M$ embedded in a sphere, and with exactly two $M$-fragments, $F$ and $F'$, such that $F$ and $F'$ do not strongly conflict. Then the embedding of $M$ extends to a flat spatial embedding of $G$, and with $F$ and $F'$ on the same side of the sphere containing $M$.
\label{noconflict}
 \end{theorem}
\begin{proof}
Assume the hypothesis and let $F=(a,b)$ and $F'=(v,w)$ be fragments that do not strongly conflict. We suppose the endpoints of $F$ and $F'$ are distinct, else $G$ is \textit {apex} (becomes planar after the removal of a vertex) and is well known to have an embedding with the desired properties. Since $G$ is connected, and $M$ is maximally planar, then $M$ is also connected. We may assume it is possible to pick disjoint paths in $M$ that connect $a$ to $b$ and $v$ to $w$. If no such disjoint paths in $M$ exist, then we pick an arbitrary path connecting $v$ to $w$, call it $P_1$, and then insert an embedded (into $S^2$) topological path into $M$, connecting, $v$ to $w$, using existing edges of $M$ where possible (intersecting an edge either entirely, or at most at its endpoints) and avoiding $P_1$, and we call this (graph theoretic) path $P_2$. We obtain an embedded supergraph of $M$, call it $M'$, and still $F$ and $F'$ do not conflict: the existence of a conflict between $F$ and $F'$ in $M'$ would require the existence of two pairs of disjoint $a,b$ and $v,w$-paths (as would be present in any graph with the required $K_{4,2}$ minor). Any such pair of paths would have to include $P_2$, so two such pairs is not possible.


We may consider $M$ on a unit sphere in ${\mathbb R}^3$. Embed each $F_i$ within a small tubular neighborhood of the corresponding $P_i$ ($i=1,2$), on the same side of the sphere containing $M$ (which we will denote $S$). Any cycle that lies in $M$ can be paneled with a disk that has interior in the opposite side of the sphere than $F$ and $F'$. We organize the rest of the (long) proof into types of cycles to be paneled.\\

\noindent
\underline{Case A}. Cycles through one fragment can be paneled:\\
\underline{Subcase A1}: a cycle through one $M$-fragment that is disjoint from the endpoints of the other $M$-fragment:\\

Consider next an $a,b$-path that avoids $v$ and $w$ that is isotopic to $P_2$ in $S^2-\{v,w\}$ (with endpoints $a$ and $b$ fixed). We will denote an arbitrary such path as $P_2'$. We claim such a path with $F$ can be paneled (with a disk that is on the same side of $S$ as $F$ and $F'$) without $F'$ getting in the way. Isotope $P_2'$ to $P_2$, keeping $F$ and $F'$ fixed. From here, it is clear that $F$ and $P_2'$ form a cycle that can be paneled, and the paneling would remain while isotoping $P_2'$ back to $P_2$; that is to say $F$ and $P_2'$ could be paneled in the original (equivalent) embedding. Similarly, every cycle $C$ that includes $F'$, avoids $a$ and $b$ with $C-F'$ isotopic to $P_1$ in $S^2-\{a,b\}$ can be paneled (with a disk on the same side of the $S$ as $F$ and $F'$) without $F'$ getting in the way. \\

Here's another justification for the skeptical reader, of the claim in the previous paragraph, that the cycle formed by $F$ and $P_2'$ can be paneled. First take a homeomorphism between $S^2-\{v,w\}$ and ${\mathbb R}^2 - \{(0,0)\}$, via stereographic projection, using $v$ as the north pole. Extend the homeomorphism so that the image of $S^2-\{v,w\}$ lies in the xy-plane of ${\mathbb R}^3$ and $F$ and $F'$ map to arcs in ${\mathbb R}^3$ that lie in $z \geq 0$ and project to $P_1$ (respectively $P_2$) in the xy-plane. Now we deform, using ambient isotopy, $F$ so that it lies in the plane $z=1$, except above $a$ and $b$, it consists of vertical segments up to $z=1$. Similarly, deform $F'$ so that it lies in the plane $z=2$, except it has vertical segments connecting up to $z=2$ at $v$ and at $w$. We apply the isotopy between $P_2$ and $P_2'$ to $F \cap (z=1)$. The result is a copy of $P_2'$ in $z=1$, which is a deformation of $F$ that could be extended to an ambient isotopy of the embedded graph, keeping all edges but $F$ fixed. This translation of $P_2'$ in $z=1$, together with the vertical segments at $a$ and $b$ and $P_2'$ in the plane $z=0$, form a cycle that can clearly be paneled with a disk that lies on the same side of $S$ as $F$. Thus the cycle formed by $F$ and $P_2'$ in the original embedding can be paneled with a disk that lies on the same side of $S$ as $F$.\\

For now, we will assume that $M$ is 2-connected. We will consider $M$ having a cut vertex later in the proof.
Now consider a cycle $C$ that contains $F$, with $P_3=C-F$ a path in $M-\{v,w\}$ connecting $a$ to $b$ that is not isotopic to $P_2$ in $S^2-\{v,w\}$. Suppose first that $P_3$ and $P_2$ meet only at $a$ and $b$. Then $P_2 \cup P_3$ separates $v$ from $w$. Because $F$ and $F'$ do not strongly conflict, then $v$ and $w$ must lie in different $P_3$-fragments. In particular, there are vertices on $P_3$, call them $x$ and $y$, such that they lie in different components of $P_3 - P_1$ and they lie in the same face of $S^2-M$. To see this, contract the $C-$fragment containing $v$ down to $v$ and the $C-$fragment containing $w$ down to $w$, deforming the embedding of $M$. Thicken $P_3$ into a topological disk, bounded by two embedded paths connecting $a$ and $b$ (sides), keeping $a$ and $b$ fixed but maintaining a planar embedding--orienting $P_3$ from $a$ to $b$ determines two sides of $P_3$, we maintain the planar embedding so that edges connecting to $P_3$ on the same side connect to the same side of the disk. There are two points on the boundary of this new rectangle, such that the two points cut the disk boundary into two components, with edges between $v$ and $C$ in one component, and edges and between $w$ and $C$ in the other. If we take these cutting vertices to be the endpoints of the minimal path in $P_3$ connecting neighbors of $v$, these two cutting vertices ($x$ and $y$ and they are not equal by the assumption of 2-connectivity) lie in a common face of $S^2-M$. Moreover, if $C=(v_1, v_2, ...,v_n)$, the cycle formed by $(v_1,..., x=v_i,y=v_k,...v_n)$ is isotopic to $P_2$ in $S^2-\{v,w\}$ (where the edge $(x,y)$ may not exist, but it can be drawn into the face and it is isotopic in $S^2-\{a,b\}$, with $x$ and $y$ fixed, to two xy-paths in $M$,  or a ``pushoff" of a path in $M$ that passes through $a$ or $b$. We will define pushoff a bit later in the proof). We can thus panel $C$ with a disk that joins to $F$, the subpath of $P_3$ from $a$ to $x$ and the subpath of $P_3$ from $b$ to $y$ on one side of the sphere, cutting through the sphere along the path from $x$ to $y$, and then joining the subpath of $P_3$ from $x$ to $y$, and thus the paneling disk avoids $F'$. See Figure \ref{proof}.

\begin{figure}[h]
\vskip -1.4in
\hskip .05in \includegraphics[scale=.53]{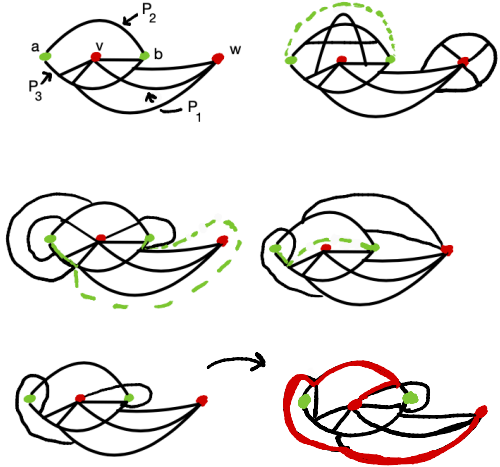}
\vskip -1.6in
\caption{Paneling around $F'$, when $P_2$ and $P_3$ meet only at their endpoints and $P_2 \cup P_3$ separates $v$ from $w$. The green path is isotopic to $P_2$ in $S^2-\{v,w\}$. The vertices $x$ and $y$ bound the dashed part of the green path. The last connected graph pictured has a strong conflict, and no such green path can be found. Highlighted in red is a $v,w$-path that, with $P_1, P_2$ and $P_3$ form a $K_{4,2}$ that establishes the strong conflict.}
\label{proof}
\end{figure}



Now suppose $P_2$ and $P_3$ meet at $a$ and at $b$ and in more vertices (and possibly edges). Then we can break up $P_3$ into a sequence of $P_2$-paths and subpaths of $P_2$. Consider one such $P_2$-path, which we will denote $P'$. Suppose $P'$ meets $P_2$ at vertices $c$ and $d$. The path $P'$ together with the subpath $cP_2d$ (which begins at $c$ and ends at $d$) forms a cycle $C'$. If $C'$ does not separate $v$ and $w$, then we do no further action for now. Otherwise, we call $C'$ a \textit{separating cycle}. We argue by induction on the number of separating cycles in $C$ that $C$ can be paneled (the previous paragraph was a sort of base case. We will discuss the other base case, when are no separating cycles, later). Similarly to the argument given in the previous paragraph, by nonconflicting of $F$ and $F'$ $v$ and $w$ cannot be in the same $P_3$-fragment. If the vertices of attachment of the fragment of $w$ to $P_3$ are contained in $P'$, or the vertices of attachment of the fragment of $v$ are contained in $P'$, it follows that there is a topological path $P''$ in $S^2-M$ that connects a point of $P'$ (call it $x$) to another point of $P'$ (call it $y$), and such that $C''=cP_2dP_3xP''yP_3c$ is a cycle that does not separate $v$ and $w$. We call $P''$ a \textit{shortcut path}. Moreover, if we can panel $C''$, then we can also panel $C=F \cup P_3$. See Figure \ref{shortcut} for a case where $P_2$ and $P_3$ meet in one vertex that is neither $a$ nor $b$. 
Otherwise, by $v$ and $w$ not in the same $P_3$-fragment, the $P_3$-fragment of $w$ cannot attach to $cP_2d$ and so the cutting vertices for the attachment of the fragment of $w$ are $x$ and $y$, where $xP_3y$ contains $cP_3d$ as a subpath. We may form a short cut path (call it $P''$ again) using $x$ and $y$ so that $P_3'=aP_3xP''yP_3b$ has fewer separating cycles that $P_3$ did, and is isotopic (in $S^2-\{a,b\}$) to a path in $M$, which has fewer separating cycles than $P_3$ and can be paneled by assumption, which implies $P_3'$ can be paneled, which implies $P_3$ can be paneled (see Figure \ref{next}). As part of this process, we may take a maximal shortcut path--one that is maximal with respect to cutting off vertices of $P_3$.


\begin{figure}[h]
\hskip 1in \includegraphics[scale=.5]{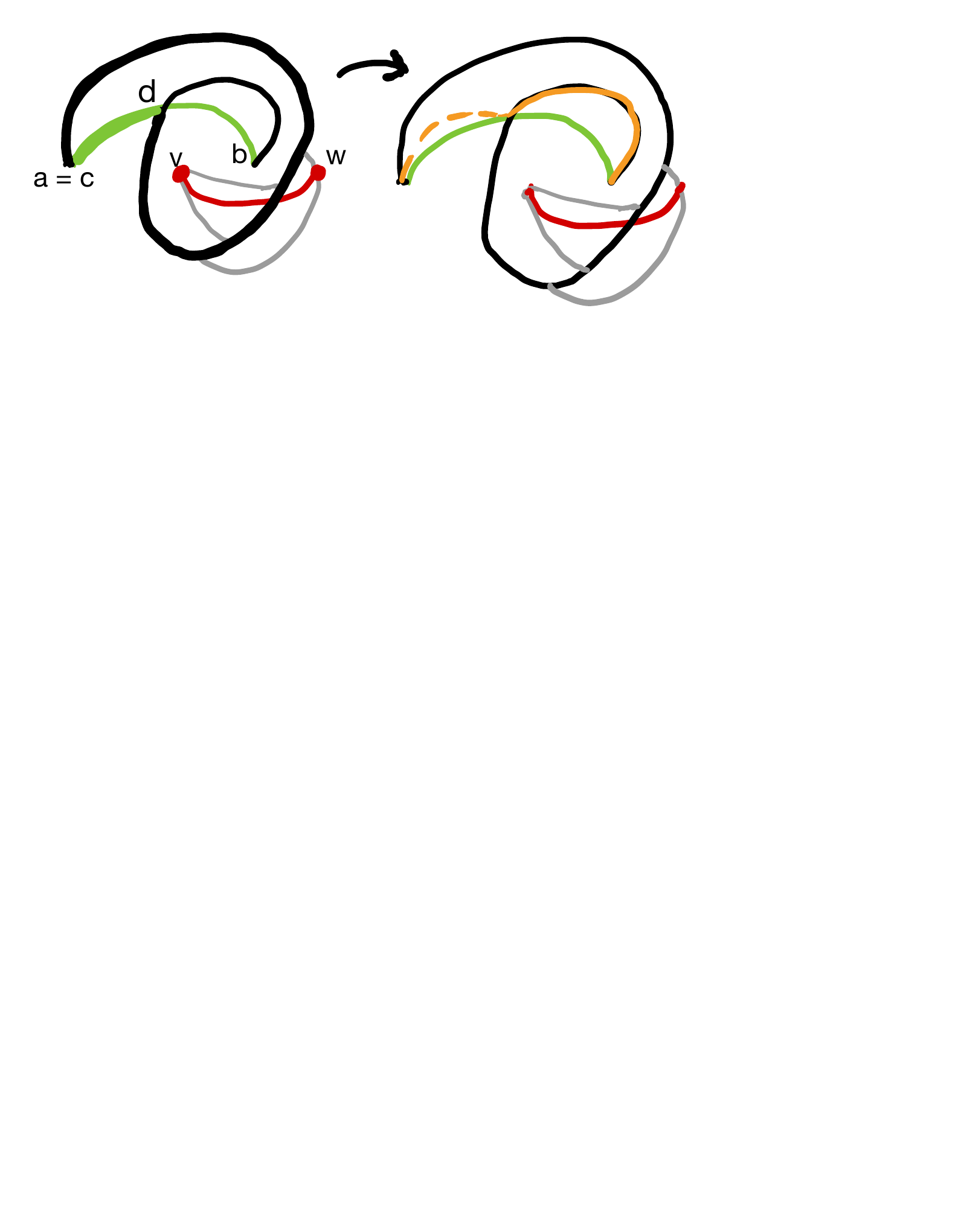}
\vskip -3in
\caption{The path $P_3$ is in black. The $P_3$-path that, with $P_2$ separates from $v$ to $w$ is in thickened and contains the vertices $c$ and $d$. On the right, the dashed path is a shortcut path, which, together with the solid orange path and $F$, forms $C''$.}
\label{shortcut}
\end{figure}


\begin{figure}[h]
\hskip .8in\includegraphics[scale=.45]{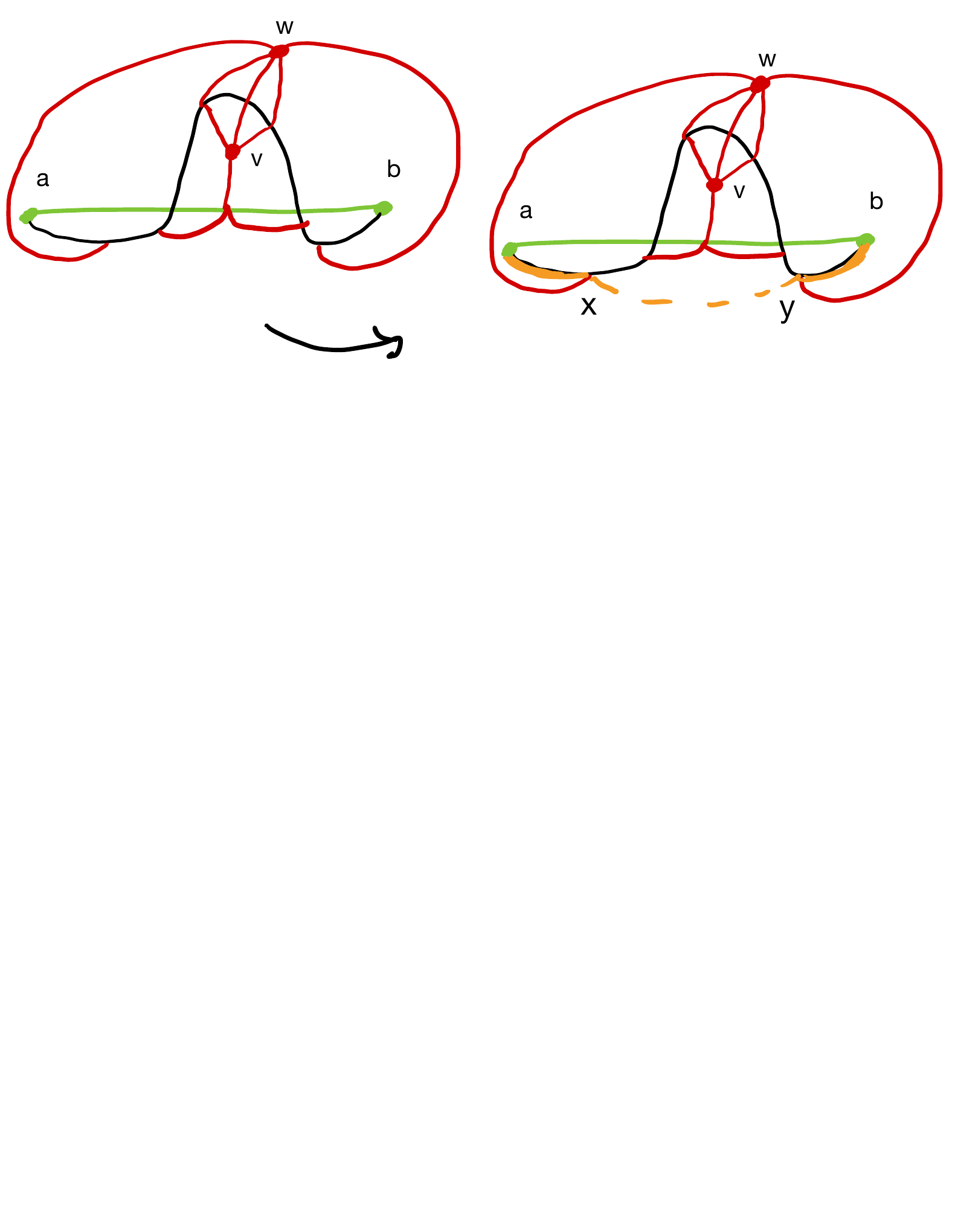}
\vskip -2.6in
\caption{A shortcut path (dashed orange) with endpoints not on the separating cycle.}
\label{next}
\end{figure}


Moreover, each shortcut path, together with a unique subpath of $P_3$ with the same endpoints, determines a cycle (a \textit{shortcut cycle}) in $M$ that can be paneled on the other side of the sphere from $F$ and $F'$. The collection of all shortcut cycles is pairwise disjoint, or intersect in a single vertex. By reasoning similar to the proof of B\"ohme's Lemma \cite{bohme}, we can simultaneously panel all shortcut cycles with disks on the other side of the sphere than $F$ and $F'$, and we can panel $F \cup P_3'$ with a disk $D$ on the same side as $F$. We can then glue $D$ to all of the disks attached to shortcut cycles to obtain a paneling for $P_3 \cup F$.

To complete the induction, we consider an $ab$-path $P_3$ that can be broken into a sequence of $P_2$-paths and subpaths of $P_2$ with no separating cycles and show that the cycle formed by such a path with $F$ can be paneled. A problematic case would be, for example, when $P_3$ intersects $P_1$. We work in ${\mathbb R}^3$ for the rest of the paragraph, with $M$ in the plane $z=0$ and $P_2$ a straight line seqment, but this is equivalent to our situation. Form a small neighborhood (in the plane $z=0$) around $P_2$ that is disjoint from $P_1$. For each $P_2$-path, we deform it into the neighborhood, taking $P_2$ paths which can be deformed (keeping $F'$ fixed) without part of $P_3$ getting in the way first. We obtain a curve that is isotopic to $P_3$, which can be paneled if and only if $P_3$ can be paneled, so we call this new path $P_3$ as well. Deform $F$ so that it has a point $R$ above the midpoint of $P_3$, and so that it consists of segments from that point to $a$ and to $b$. The cone of the path to this point $R$ gives a paneling disk. It will not intersect $F'$, as $F'$ is not above this neighborhood. This completes our induction.\\

 \noindent
 \underline{Case A2:} a cycle through one fragment that is not disjoint from one or both endpoints of the other fragment:\\
 
We next consider any cycle $C=F \cup P$, where $P$ includes $v$ or $w$. We use the same shortcutting process. The resulting short cut version of the path meets $P_1$ only possibly at $v$ or at $w$. At such a point of intersection, say $v$, we first consider a small neighborhood around $v$ (so small that it meets the path only in the two edges incident to $v$), and then taking the intersection of the boundary of that neighborhood with the path--to get two points, $x$ and $y$, and replacing the part of the path that is inside the neighborhood, with an arc of the intersection that is disjoint from $P_1$ and connects $x$ and $y$. Repeating this for $w$, if necessary, we obtain what the call the \textit {pushoff path}. As above, this pushoff path can be isotoped to a neighborhood of $P_2$. It follows, as above, that it can be paneled with a disk that lies on the same side of the $S^2$ and $F$. From there, it follows, by pushing the pushoff path back to $P$, that $F$ with the short cut path can be paneled with a disk that lies on the same side of the sphere. Finally, we have that $F \cup P$ can be can be paneled.\\

 \noindent
 \underline{Case B:} a cycle through both fragments:\\
 
 Finally, suppose that $F$ and $F'$ lie in a common cycle, $C$. The cycle consists of $F$, $F'$ and $P$ and $P'$, which are paths in $M$. The cycle $C$ cannot be part of a nonsplit link, as every possible disjoint cycle can be paneled, by our above arguments. As $M$ is connected, there is a path from $P$ to $P'$ and we can use this path to cut $C$ into a theta-curve, two cycles of which can be paneled (they do not contain both fragments). Since their intersection is connected, we can panel both simultaneously (see B\"ohme's Lemma). It follows that $C$ is unknotted. It remains to show that $C$ can be paneled.
 
 To show $C$ can be paneled, we use some known results from the theory of tangles. From \cite{fl} a \textit{tangle} is a three-dimensional ball $B$, containing two disjoint arcs together with a finite number (for us zero) of disjoint simple closed curves that are all represented by $t$, such that the intersection of $t$ with the boundary of $B$ is precisely the set of endpoints of the two arcs. Two tangles $(B, t_1)$ and $(B, t_2)$ are equivalent if there is an orientation preserving homeomorphism $h:(B, t_1) \rightarrow (B, t_2)$ that is the identity on the boundary of $B$. 
 
 
 We say a tangle $(B,t)$ is \textit{rational} if $t$ can be deformed by an ambient isotopy in $B$ to the trivial tangle, where the endpoints of $t$ remain in the boundary of $B$ throughout the isotopy. Conway \cite{conway} defined the \textit{fraction} of a rational tangle
to be a rational number or $\infty$ and proved that two rational tangles are equivalent if and only if they have the same fraction (See Figure \ref{triv} for drawings of the $0$ and $\infty$ tangles).
 
 From a standard projection of a (rational) tangle, there are two common ways to close the tangle into a knot (or link): the numerator and the denominator (see Figure \ref{closure}). 
 
 \begin{figure}
 \hskip .8in \includegraphics[scale=.5]{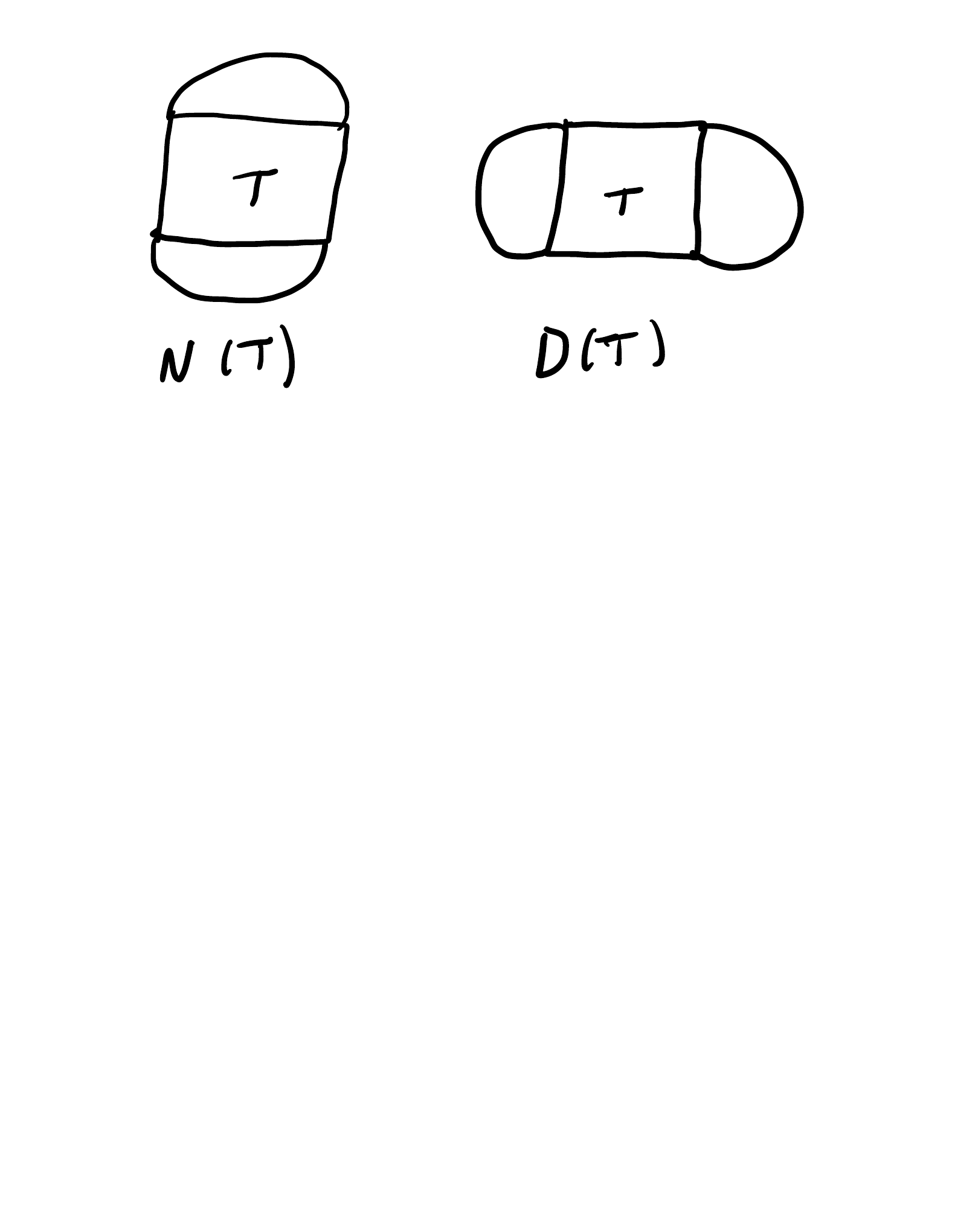}
 \vskip -2.75in
 \caption{The numerator and denominator closures of a tangle $T$, $N(T)$ and $D(T)$.}
 \label{closure}
 \end{figure}

 As suggested by Colin Adams, we will use Schubert's classification of rational knots (\cite{schu} see also Kauffman and Lambroupoulu \cite{kaufL}): Suppose that rational tangles with fractions $p/q$ and $p'/q'$ are given ($p$ and $q$ are relatively prime. Similarly for $p'$ and $q'$). If $K(p/q)$ and $K(p'/q')$
 denote the corresponding rational knots obtained by taking numerator
closures of these tangles, then $K(p/q)$ and $K(p'/q')$ are equivalent if and only if:\\
1. $p = p'$\\
2. either $q \equiv q'$  $mod$ $p$ or $qq' \equiv 1$ $mod $ $p.$\\

 The fragments $F$ and $F'$, being placed near the disjoint paths $P_1$ and $P_2$ can considered as forming a trivial tangle. We can further consider the knot $C$ through $F$ and $F'$ as being formed as the closure of a related tangle. We deform $P_1$ and $P_2$ in the boundary of the ball until they look like the trivial tangle 0 (see Figure \ref{triv}). Since the endpoints of $F$ and $F'$ stayed on the boundary, extending the deformation to the interior of the ball results in $F$ and $F'$ forming a rational tangle with $C$ being the numerator closure of the tangle. Since $C$ was itself unknotted, by Shubert's theorem, it follows that $C=K(p/q)$ is equivalent to $K(1/1)$ (see Figure \ref{one}) or $C$ was the numerator closure of the tangle $\infty$. If $C=K(\infty)$ (see Figure \ref{triv}), then by direct inspection, one can see $C$ can be paneled by a disk on the same side of $S$ as $G$. If $C$ has a rational number fraction $p/q$, it follows that $p=1$ and $q$ is an integer (nonzero).
 
 \begin{figure}
 \hskip 1.2in \includegraphics[scale=.5]{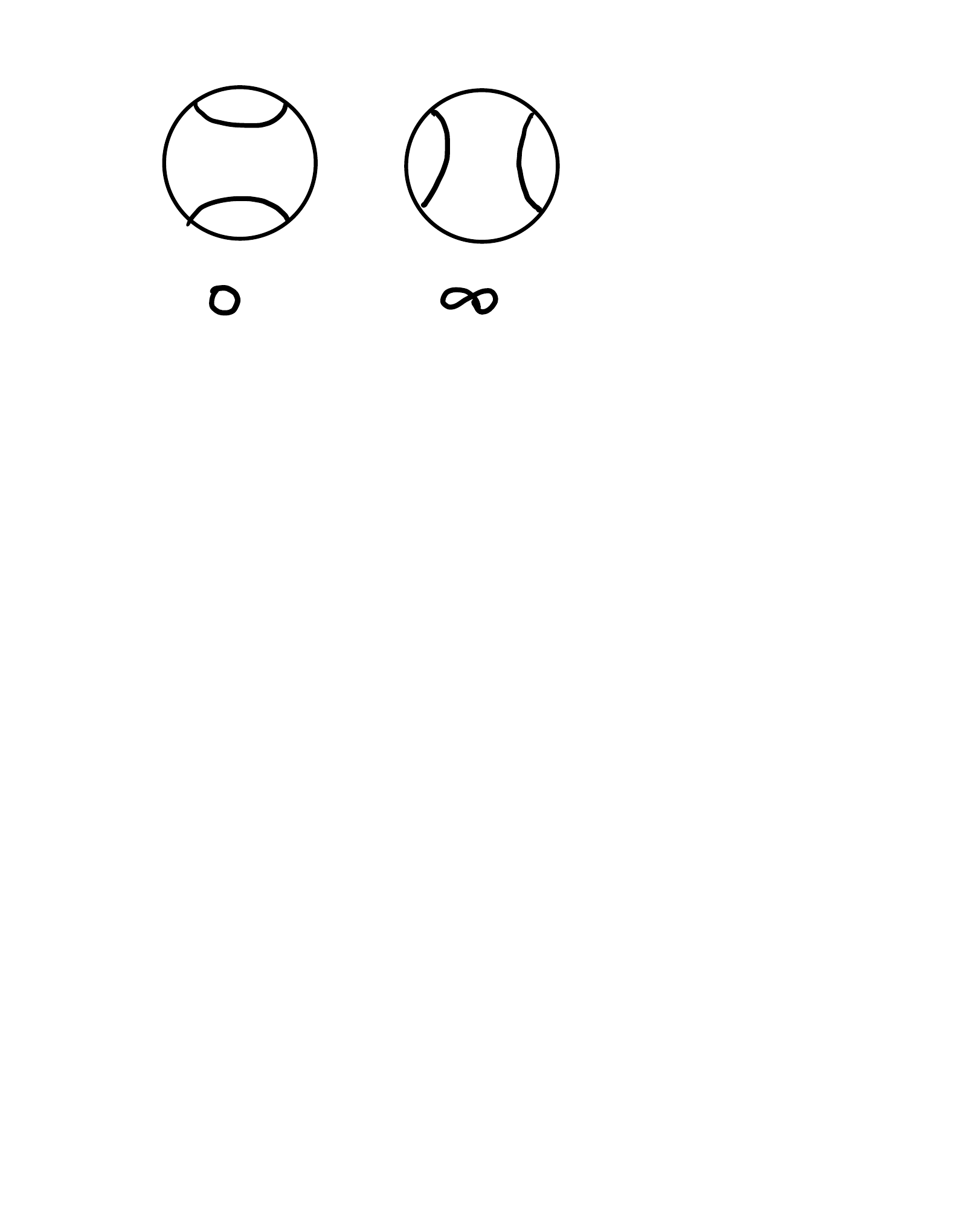}
 \vskip -2.95in
 \caption{The 0 and $\infty$ tangles.}
 \label{triv}
 \end{figure}

 To finish the argument, we need to examine the denominator closure of $T$. This knot (link) is present in the embedding of $M \cup F \cup F'$, as the paths $P_1$ and $P_2$ are disjoint and connect $a$ to $b$ and $v$ to $w$. It follows that a subpath of $P_1$ and a subpath of $P_2$ connect $P$ and $P'$. The subpaths, taken with $P$ and $P'$ form a cycle with distinct vertices $a', w', b', v'$ with either $a'=a$ or $a'$ connected to $a$ by a subpath of $P$ (similarly, either $b'=b$ or $b'$ connected to $b$ by a subpath of $P'$, either $v=v'$ or they are connected by a subpath of $P$ and finally, either $w=w'$ or they are connected by a subpath of $P'$). On the 4-cycle, there are two possible orderings of vertices, going, say, clockwise: $a',w',b', v'$ or $a', b', w', v'$. In the first case, the subpath of $P$ from $a'$ to $v'$ and the subpath of $P'$ from $w'$ to $b'$ forms the denominator closure of the tangle formed by $F$ and $F'$. Note that $D(T)=N(T^r)$ (\cite{kaufL}), where $T^r$ is the rotation of $T$, 90 degrees counterclockwise. As $N(T^r)$ is trivial, this implies that the fraction associated to $T^r$ is of the form $1/n$, where $n \in {\mathbb Z} -\{0\}$, or $\infty$. If $T^r=K(\infty)$, then $T=K(0)$, but $T$ is not a 2-component link, so $T^t=K(1/n)$. Since $T^r=-1/T$ \cite{kaufL}, the fraction associated to $T^r$ is $-f(T)$, which, (if $p/q$ is rational), combining equations gives: $1/n=-q/1$, where $n$ and $q$ are integers. But then $-nq=1$, so either $n=1, q=-1$ or $q=1, n=-1$. In either case, $T$ has fraction $-1$ or $1$ which one can check that $C$ can be paneled with a disk that meets $M$ only along $P_1$ and $P_2$ (see Figure \ref{one}). That is, $C$ is paneled.

 \begin{figure}
\hskip 1.4in \includegraphics[scale=.4]{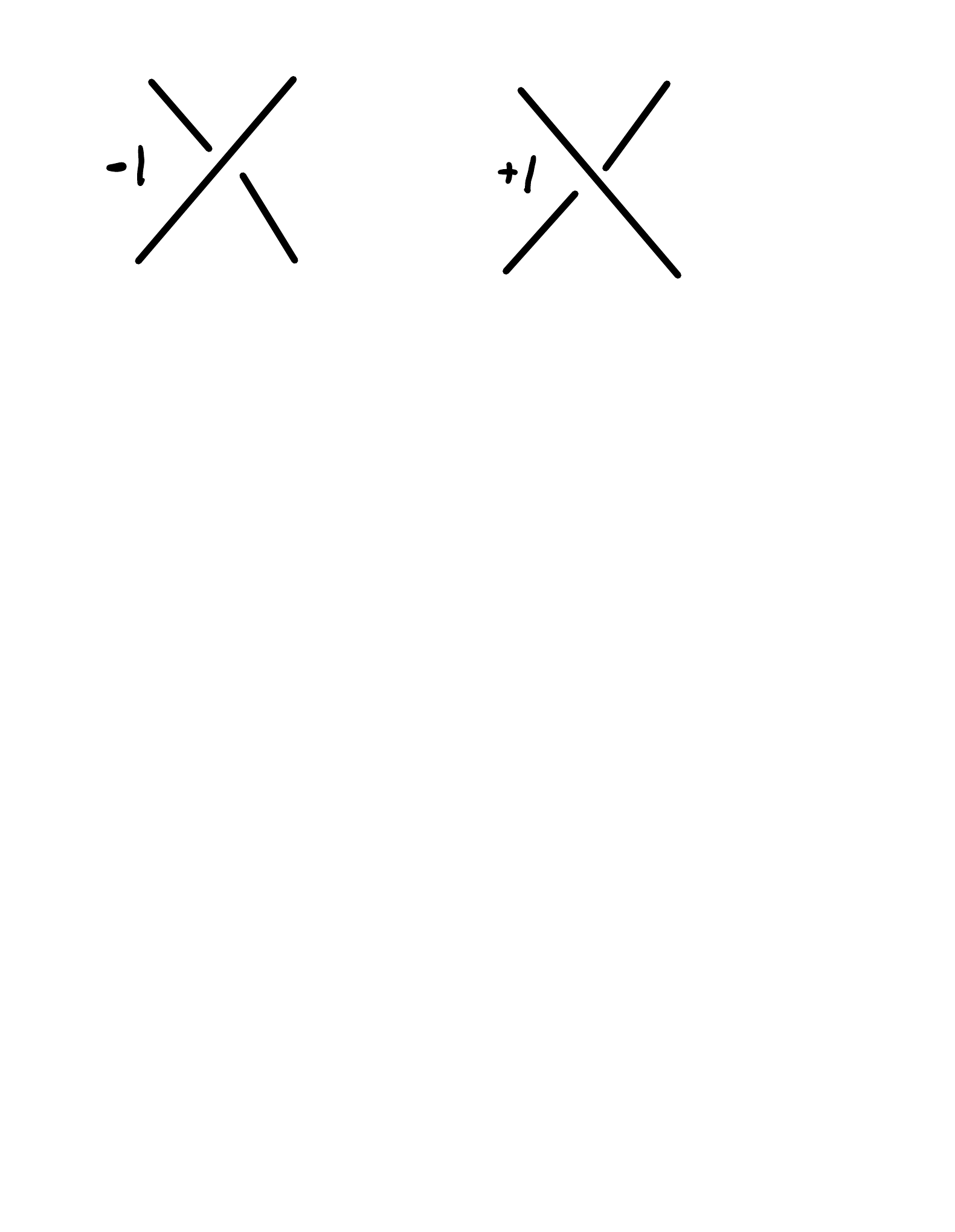}
 \vskip -2.4in
 \caption{The tangles -1 and +1.}
 \label{one}
 \end{figure}
 There remains the case with the vertices appearing $(a',b',w',v')$. In this case, $D(T)$ is a split link of two components, so its fraction equals to $0$. We have $D(T)=N(T^r)=-N(1/T)=0$, which implies that $N(T)$ is infinity, which means $T$ is not equivalent to any other $p/q$ and $N(T)$ must be panelable, as we discussed above.
 

As $C$ was an arbitrary cycle through both $F$ and $F'$, then it follows that the embedding of $G=M\cup F \cup F'$ is flat.\\

\underline{Final Loose End:} When $M$ has a cut vertex:\\

Finally, we consider what happens when $M$ has a cut vertex. We break this up into a few cases. First, if $P_1$ and $P_2$ lie in a common block, then there is a flat embedding of $M\cup F \cup F'$ by the arguments given above (the other blocks only contribute cycles that lie on $S$ which can clearly be paneled and the other blocks easily get out of the way of panels of cycles in the block that contains $P_1$ and $P_2$.) Otherwise, suppose next there is a cut vertex that does not lie in the interior of $P_1$ nor in the interior of $P_2$ (if it lies at an endpoint of $P_i$, it cuts the other endpoint off from the other $P_j$ , $i,j=1,2, i \neq j$, else we have the previous case). Then every cycle through $F'$ corresponds to a $v,w$-path, and every $v,w$-path is isotopic to $P_1$ in $S^2-\{a,b\}$. Similarly, every $a,b$-path is isotopic to $P_2$ in $S^2-\{v,w\}$ and every cycle through $F$ can be paneled, as discussed earlier.  There may be a $v,w$-path that passes through $a$ (or $b$), but its push off path is isotopic to $P_1$, and thus the cycle it forms with $F'$ can be paneled. Similarly for $a,b$-paths through $v$ (or $w$). Because of the cut vertex, there is no cycle that contains $F$ and $F'$.

Finally, suppose the cut vertex lies in an interior point of $P_1$ (or $P_2$--for now suppose $P_2$ has no cut vertex). The cut vertex cuts $M$ into at least two components, one containing $P_2$. Either $v$ or $w$ lies in a different component than $P_2$. Suppose $w$ lies in component $H$ that does not contain $P_2$. Consider the embedded subgraph of $M$ formed by removing all vertices and edges of $H$, except those on $P_1$. By assumption, $P_1$ and $P_2$ lie in a block in this embedded subgraph, which can be embedded flatly with $F$ and $F'$. Putting the $H$ back in the embedding creates some cycles that are in $S^2$ and so panelable, as well as some new $v,w$-paths (cycles through $F$), but because of the cut vertex, such paths will be isotopic to $P_1$ in $S^2-\{a,b\}$. As before, such cycles can be paneled. Similar reasoning applies if $P_2$ contains a cut vertex of $M$. This concludes our proof.
\end{proof}

We can generalize the previous result to more than 2 fragments that do not pair-wise conflict:
\begin{corollary}

Let $G$ be a non-planar graph with (maximally) planar subgraph $M$ embedded on a sphere, and with exactly $n$ $M$-fragments, $F_1, F_2, ...,F_n$ such that $F_i$ and $F_j$ do not strongly conflict for $i\neq j$, with $n \geq 2$. Then the embedding of $M$ extends to a spatial embedding of $G$, and with all of the $F_i$ on the same side of the sphere containing $M$, and such that the embedding restricted to $M$ union any two fragments is flat.
\label{forn}
 \end{corollary}
 
 \begin{proof}
Note that there is a unique, up to ambient isotopy, way to embed an unknotted segment inside a sphere, with endpoints on the sphere. We further note that for the $n=2$ case (previous proof), we can embed the first fragment, $F$, and consider it fixed, then embed the second fragment, $F'$, and then adjust--pushing $F'$ through $F$ as necessary, so that the embedding of $M \cup F \cup F'$ is linkless. Now suppose $n>2$ and the result is true for $n-1$ such fragments.  Consider the first $n-1$ embedded fragments as fixed. For $F_n$, embed it in the sphere unknotted, and then adjust with $F_1$, so that $F_1 \cup F_n \cup M$ forms a linkless embedding. We can then adjust $F_n$ with $F_2$--pushing $F_n$ through $F_2$ as necessary, but not allowing $F_n$ and $F_1$ to pass through each other. We repeat this process up to $F_{n-1}$, until the embedding restricted to $M$ union any two fragments is linkless. 
 \end{proof}

\begin{figure}
\hskip 1.2in
\includegraphics[scale=.165]{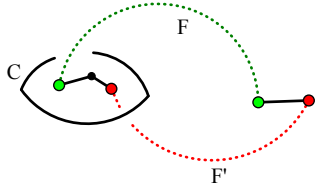}
\vskip -.2in
\caption{A nonsplit link occurs when strongly anti-conflicting $M$-fragments are placed on different sides of a sphere containing $M$.}
\label{fig:anticonflict}
\end{figure}

We introduce another definition. 
Let $G$ be a graph with $M$ a maximally planar subgraph that is embedded in $S^2$. Let $F=(v_1,v_2)$ and $F'=(w_1,w_2)$ be two fragments of $M$. We say that $F$ and $F'$ {\bf strongly anti-conflict} if there exists a cycle $C$ in $M$, with the vertices of attachment of $F$ ($F'$) lying in different components of $S^2 -C$, and there is a path (possibly trivial) in $M$, connecting $v_1$ and $w_1$ in $S^2 -C$ and a path (possibly trivial) in $M$ connecting $v_2$ and $w_2$ in $S^2 - C$. It is clear from Figure \ref{fig:anticonflict} that if two strongly anti-conflicting $M$-fragments are placed on opposite sides of a sphere containing $M$, then the resulting embedding contains a pair of nonsplittably linked cycles. We call this the \textit{Anti-conflicting Observation}.

\begin{prop}

Let $G$ be a connected non-planar graph with maximally planar subgraph $M$ embedded in a sphere, and with exactly two $M$-fragments, $F$ and $F'$, such that $F$ and $F'$ do not strongly anti-conflict. Then the embedding of $M$ extends to a flat spatial embedding of $G$, and with $F$ and $F'$ on opposite sides of the sphere containing $M$.
\label{noconflict2}
 \end{prop}
 
 \begin{proof}
Assume the hypotheses, and place $F$ and $F'$ on opposite sides of $S$, the sphere containing $M$, each, together with a path in $M$ connecting its endpoints, bounding a disk with interior on the same side of $S$ as $F$ ($F'$). Any cycle through $F$ that does not contain $F'$ can be paneled, just as any cycle in an apex graph embedding can be paneled. This is seen by removing $F'$: the resulting embedding is apex. Placing $F'$ back will not puncture the paneling disk, as $F'$ is on the other side of the sphere. Similarly, any cycle through $F'$ that does not contain $F$ can be paneled. \\

\noindent
\underline{A cycle that lies in $M$ can be paneled:}\\

If a cycle lies entirely in $M$, we claim that it can be paneled. Let $C$ be an arbitrary such cycle. Denote $F=(v_1,v_2)$ and $F'=(w_1,w_2)$. If both endpoints of $F$ or of $F'$ lie in the closure of the same component of $S^2-C$, then the result follows, as the cycle can be paneled by a disk that lies on the same side of the sphere as the fragment with both endpoints in the closure of the same component of $S^2-C$: connect the endpoints of the fragment by an embedded topological path (path whose image is homeomorphic to a unit interval) with interior in that component of $S^2-C$, then deform the newly made path to the other side of $S^2$, and then deform the fragment so that it lies in a small neighborhood of this path, out of the way of a paneling disk on the original side of $S^2$.

 We now assume both fragments have endpoints on different sides of $S^2-C$. We may take $v_i$ on the same side as $w_i$. As $F$ and $F'$ do not strongly anti-conflict, without loss of generality, we may assume that every path in $M$ from $v_1$ to $w_1$, must pass through $C$. That is to say $v_1$ and $w_1$ lie in different $C$-fragments in $M$.

In what follows, when we refer to a $C$-fragment, we mean $C$-fragment in $M$ (ignoring $F$ and $F'$). If the $C$-fragment of $v_1$ connects to $C$ in exactly one vertex, or if it is not connected to $C$ at all, then the spatial embedding of $M\cup F \cup F'$ is ambient isotopic to an embedding for which $v_1$ and $v_2$ are on the same component of $S^2-C$ (flip the $C$-fragment of $v_1$ over to the other side of $C$ on $S^2$), which we have already discussed. We now assume that the $C$-fragment of $v_1$ and the $C$-fragment of $w_1$ connect to $C$ in at least two vertices.

Similar to the proof of Theorem \ref{noconflict}, the $C$-fragment of $v_1$ and the $C$-fragment of $w_1$ help us determine our paneling disk. Note that  the $C$-fragment of $v_1$ and the $C$-fragment of $w_1$ cut $C$ up into four paths (two possibly trivial): denote by $C_1$ the minimal path contained in $C$ that contains all vertices attached to the $C$-fragment of $v_1$, such that $C_1$ does not contain any vertices attached to the $C$-fragment of $w_1$. Then we denote by $C_2$ the minimal path contained in $C$ that contains all vertices attached to the $C$-fragment of $w_1$, and no vertices attached to the $C-$fragment of $v_1$. We denote $C_3$ and $C_4$ the possibly trivial paths that connect an endpoint of $C_1$ to an endpoint of $C_2$, with $C_i \cap C_j$ at most one vertex, for $i=3,4$, $j=1,2$. Since $v_1$ and $w_1$ lie in different $C$-fragments, the endpoints of $C_1$, and similarly the endpoints of $C_2$, must lie in the same face of $S^2-M$. Form an embedded topological path in the face that connects the endpoints of $C_i$, and call the path $P_i$, $i=1,2$. We have that $P_i \cup C_i$ is a circle, and we can panel it with a disk, $D_i$ that is on the opposite side of $M$ as $F$ for $i=1$ and on the opposite side of $M$ as $F'$ for $i=2$. The circle $C'$ formed by $P_1, C_3, P_2, C_4$ can also be paneled with $D_3$, as there are no endpoint of $F$ nor $F'$ on one side of $S^2- C'$. Finally, we can form a paneling disk for $C$ with $D_1 \cup D_2 \cup D_3$, see Figure \ref{panel}. Thus all cycles in $M$ can be paneled.\\

\begin{figure}
\hskip 1in \includegraphics[scale=.65]{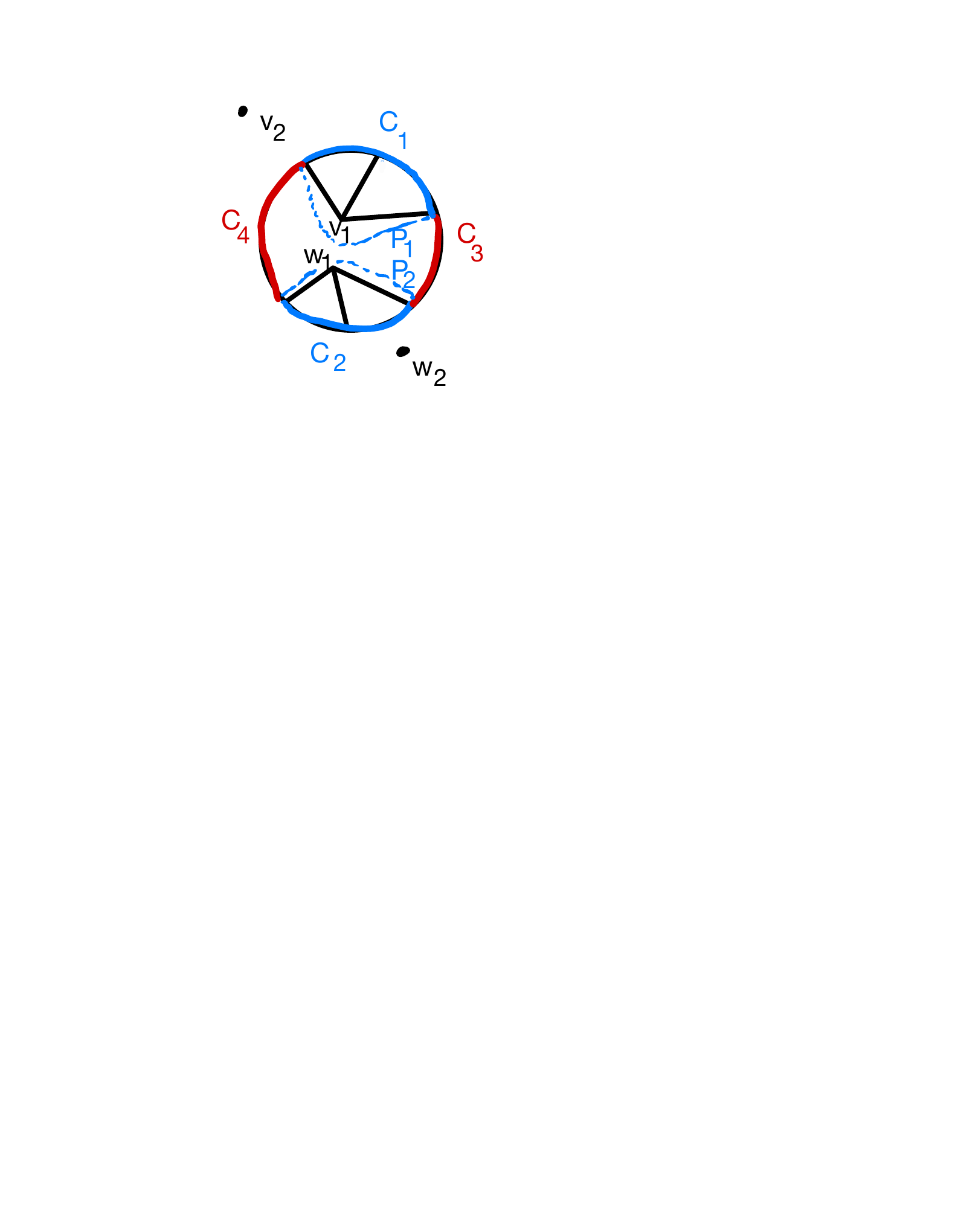}
\vskip-3.6in
\caption{Paneling a cycle in $M$.}
\label{panel}
\end{figure}

\noindent
\underline{A cycle that contains both $F$ and $F'$ can be paneled:}\\

Next, if a cycle contains $F$ and $F'$, then it must be of the form $C=P_1, F, P_2, F'$, where $P_1$ and $P_2$ are paths in $M$ joining an endpoint of $F$ to an appropriate endpoint of $F'$. In terms of cycles that possibly form a nonsplit link with $C$, the only cycles that are disjoint from $C$ lie entirely on $S$ (in $M$). However, we showed earlier in this proof that every such cycle can be paneled. Thus $C$ does not link with any other cycle. 

We claim further that $C$ can be paneled. Since, by definition of not strongly anti-conflicting, there is no cycle in $M$ separating $P_1$ from $P_2$, we will show after this paragraph that there is an embedded topological path, $\gamma$, with interior in $S-M$ from $P_1$ to $P_2$. We panel the cycle formed by $F$, $\gamma$, and part of $P_1$ and $P_2$. We similarly panel the cycle formed on the other side of $M$ by $F'$, $\gamma$, and the rest of $P_1$ and $P_2$. We can glue these paneling disks together along $\gamma$ to get a paneling disk for $C$. This concludes our proof, once we check the path claim.\\

To check our claim that there is an embedded topological path, $\gamma$, with interior in $S-M$ from $P_1$ to $P_2$; take the union of all closures of path components of $S-M$ [with their boundaries], that intersect $P_2$, and call this union $U$. We show that if $U$ does not intersect $P_1$, then there must have been a cycle separating $P_1$ and $P_2$, which contradicts $F$ and $F'$ not strongly anti-conflicting. Then $U$ intersecting $P_1$ implies the existence of such a path $\gamma$. We provide more details in what follows.

 Consider the boundary of $U$ in $S$. It is a subset of $M$, for if $p \in S$ with $p \notin M$, then there is a small ball around $p$ that is disjoint from $M$, which must be contained in either the interior of $U$ or the interior of $U^c$. Then $p \notin Bd(U)$, and so $Bd(U)\subseteq M$.


 Assume, for the sake of contradiction, that $P_1$ is disjoint from  $U$. It follows that $P_1$ lies in a component of $S-U$, which we will denote as $V$. This component must be open as it is a component of an open set in a locally connected space.

Since $V$ is open and non-empty and a proper subset of $S$, $V$ is not closed, and so $Bd(V) \neq \emptyset$. We have $Bd(V) \subseteq M$ by a similar argument as was used to show $Bd(U) \subseteq M$. If $Bd(V)$ does not have any edges, then for some vertex $p \in Bd(V)$, there is a small neighborhood of $p$ that intersects $M$ only in $p$ and edges adjacent to $p$. Each sector determined by the edges incident to $p$ and the neighborhood, must either be in $U$ or in $S-U$. As $p \in Bd(V)$, some sector must be a subset of $V$. Start with that sector and go clockwise. Eventually, as $p \in Bd(V)$, we must get a sector contained in $V$, followed by a sector contained in $V^c$, but then the edge between these sectors should be included in $Bd(V)$. This is a contradiction, so $Bd(V)$ must contain an edge. 

Pick one such edge. Orient across the edge so that $V$ is on the left and $V^c$ is on the right. At the terminal vertex of this edge, we have a similar situation, where the first edge clock-wise from our first edge has $V$ on the left and $V^c$ on the right (else the $V$ and $V^c$ regions would intersect). We continue this way, forming a path, until we have a repeated vertex at the end of a new edge. The repeated vertex must be the first vertex of the first edge, or else again $V$ and $V^c$ regions would intersect. Thus we have formed a cycle that separates $P_1$ from $P_2$, a contradiction that $F$ and $F'$ do not anti-conflict.

It follows that $P_1$ and $P_2$ both intersect the boundary of a common path component of $S-M$. Thus, there is an embedded topological path in $S$ from $P_1$ to $P_2$ that meets $M$ only in its endpoints.

\end{proof}

\section{The conflict graph}

  In order to generalize the notion of conflict to spatial graphs, we must discuss the notion of signed graph. (Thanks to Garry Bowlin for first suggesting signed graphs would be the appropriate generalization here). A \emph{signed graph} $\Sigma = ( G, \sigma )$ consists of a graph $G$ and an edge signing $\sigma : E( G ) \to \{ +, - \} $. The sign of a cycle $C=e_1e_2\cdots e_{n-1} e_n$ is obtained by multiplying the signs of its constituent edges:
$$
\sigma (C) = \sigma (e_1) \sigma (e_2) \cdots \sigma (e_n)
$$
A signed graph $\Sigma$ is \emph{balanced} if all of its cycles are positive. Otherwise, it is \emph{unbalanced}. Harary \cite{harary} showed that a signed graph is balanced if and only if its vertex set can be divided into two sets (one of which may be empty), $X$ and $Y$, so that each edge between the sets is negative and each edge within the sets is positive. In this way, balanced signed graphs are a generalization of bipartite graphs.

Let $M$ be a maximally planar subgraph of $G$, embedded in a sphere.  The \textit{strong conflict graph} of $M$ is the graph whose vertices are fragments of $M$.  Two vertices in the strong conflict graph are connected by a positive edge if their corresponding fragments strongly anti-conflict and by a negative edge if they strongly conflict.

\begin{figure}[h]
\vskip -.25in
\hskip 1.1in  \includegraphics[scale=.21]{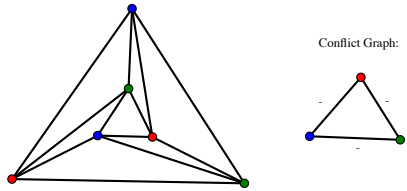}
\vskip -.2in
\caption{A maximal planar subgraph of $K_6$ and its associated (strong) conflict graph. Negative signs indicate strongly conflicting fragments. \cite{fr}}
\label{Pet1}
\end{figure}

\begin{figure}[h]
\vskip -.05in
\hskip .9in
\includegraphics [scale=.225]{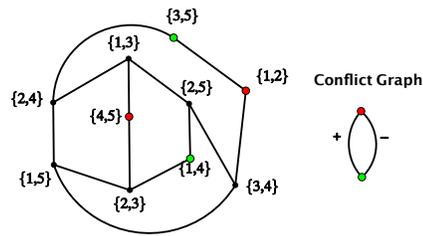}
\vskip -.3in
\caption {A maximal planar subgraph of the classic Petersen graph and its associated (strong) conflict graph. A positive edge indicates a strong anti-conflict. \cite{fr}}
\label{Pet2}
\end{figure}

\begin{figure}
\vskip -.4in
\hskip 1.2in
\includegraphics[scale=.2]{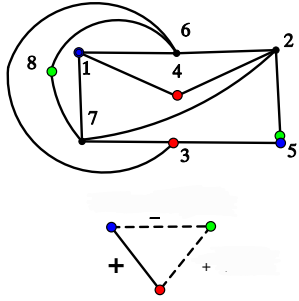}
\vskip -.3in
\caption{A maximal planar subgraph of $K_{4,4}-e$ (\cite{fr}) with an associated strong conflict graph that is balanced. See \cite{fr} for more details as to why the associated conflict graph is not balanced (dashed edges indicate implicit (anti-)conflict).} 

\label{strongbad}
\end{figure}

As the example in Figure \ref{strongbad} suggests, for $K_{4,4}-e$, two fragments of a maximally planar subgraph may not strongly conflict, but after a suitable adjustment, the fragments strongly conflict as fragments of a related planar subgraph. There are three other strong conflict graphs of $K_{4,4}-e$ that are balanced, but with a broader definition of conflict, one can show that the resulting conflict graphs are all unbalanced. For every other graph in the Petersen Family, every associated strong conflict graph is unbalanced. See Figures \ref{Pet1} and \ref{Pet2}, for example. See the following discussion and \cite {fr} (Proposition 7) for more details.

\vskip .4in




 If the strong conflict graph is balanced, we use the balanced strong conflict graph to place the fragments (unknotted with respect to the sphere) on the appropriate side of the sphere: conflicting on opposite sides of the sphere, anti-conflicting on the same side, and indifferent edges (that are neither conflicting nor anti-conflicting) arbitrarily with respect to each other. For two fragments that do not conflict, if they are on the same side of the sphere, we place the fragments so that, together with the planar subgraph, they pair-wise form a flat embedding (ignoring the other fragments: see Corollary \ref{forn}). We call such an embedding a \textit{potentially flat} embedding. There may be multiple non-equivalent potentially flat embeddings. 

We say non-strongly conflicting (respectively non-strongly anti-conflicting) $M$-fragments $F$ and $F'$ {\bf implicitly conflict (resp. implicitly anti-conflict} if, for every potentially flat embedding with $F$ and $F'$ on the same side of the sphere containing $M$ (resp. on different sides), there exists a sequence of $S^2$-embedded graphs $M=M_1, ..., M_n$ where
 each $M_i$ results from $M_{i-1}$ from a single edge contraction or removal, $M_i$ maintains the embedding into $S^2$ inherited from $M_{i-1}$ and $M_i$ accepts any entire fragments possible by ambient isotopy, with $F$ and $F'$ staying on the same side of the $S^2$ throughout, with $F$ and $F'$ strongly conflicting (strongly anti-conflicting) with respect to $M_n$.
 
 Given a nonplanar graph $G$ with planar subgraph $M$ embedded on a sphere, the \textit {conflict graph} of $M$ has the fragments of $M$ as vertices and includes negative edges between strongly or implicitly conflicting fragments, and positive edges between strongly or implicitly anti-conflicting fragments.
 
 
 For a given potentially flat embedding of $G$, there are only finitely many different possible nonequivalent $S^2$ embeddings of minors of $G$ (possible $M_i$). It follows from the definition of conflict graph as well as Proposition 2.1 and the Anti-conflicting Observation, that if the conflict graph associated to a graph $G$ is unbalanced, then $G$ does not have a (spherically) flat embedding. Conversely, if we start with a planar subgraph $M$ of a graph $G$ with flat embedding, then a conflict graph of $M$ must be balanced. We summarize:
 

\begin{theorem}
If $G$ has a flat embedding, then within a flat embedding, for an arbitrary maximal planar subgraph of $G$, embedded on a sphere, the associated conflict graph must be balanced. Equivalently, given a graph $G$ with maximal planar subgraph $M$, if for every embedding of $M$ into a sphere, the conflict graph of $M$ is unbalanced, then $G$ is has no (spherically) flat embedding.
 \end{theorem}

 It has also been shown:
 
 \begin{theorem} \cite{fr} Every maximal planar subgraph of a Petersen Family graph has every possible conflict graph unbalanced. There are 45 such maximally planar subgraphs. 
\end{theorem}



\section*{Acknowledgments}
This work started out of a 2012 REU project with students Cara Nickolaus, Justin Raimondi, Joshua Wilson and Liang Zhang. I wish to thank the SUNY Potsdam-Clarkson University REU Program and the NSA and NSF for financial support under NSA Grant H98230-11-1-0206 and NSF Grant DMS-1004531. I would like to thank Colin Adams, Garry Bowlin, Ramin Naimi and Tom Zaslavsky for helpful conversations.


{\it Joel Foisy, Mathematics Department, SUNY Potsdam, Potsdam, NY 13676, foisyjs@potsdam.edu\\


\begin{thebibliography}{999}

\bibitem{alex} Alexander, J.W., On the subdivision of 3-space by a polyhedron. \textit{Proc. Nat. Acad. Sci. U.S.A.} 10, 6-8, 1924.

\bibitem{bohme}
B\"ohme, T., \textit{On spatial representations of graphs}, Contemporary Methods in Graph Theory (R. Bodendieck, Ed.), Mannheim, Wien, Zurich, (1990), 151-167.

\bibitem{bl} 

Bonnington, C.P. and Little, C.H.C., The foundations of topological graph theory, Springer-Verlag, New York, 1995.



\bibitem{conway} Conway, J.H., An enumeration of knots and links and some of their algebraic properties, Proceedings of the conference on Computational problems in
Abstract Algebra held at Oxford in 1967, J. Leech ed., (First edition 1970),
Pergamon Press, 329-358.

\bibitem{conwaygordon}
		Conway, J.H. and Gordon, C.  \textit{Knots and links in spatial graphs}. J. Graph Theory 7 (1983), no. 4, 445-453.
		
		\bibitem{DS2} Dirac, G. A.; Schuster, S.
A theorem of Kuratowski.
Nederl. Akad. Wetensch. Proc. Ser. A. 57 = Indagationes Math. 16, (1954). 343-348.

\bibitem {DS} Dirac, G. A.; Schuster, S. Corrigendum: ``A theorem of Kuratowski,'' Nederl. Akad. Wetensch. Proc. Ser. A 64 = Indag. Math. {\bf 23} (1961) 360.
\bibitem{fl} Flapan, Erica, When topology meets chemistry. (English summary)
A topological look at molecular chirality. Outlooks. Cambridge University Press, Cambridge; Mathematical Association of America, Washington, DC, 2000.

\bibitem{fr}
Foisy, J. and Raimondi, J. \textit {Conflict graphs of maximal subgraphs of the Petersen graphs}, arXiv:2207.06251, 2022.

\bibitem{harary}
 Harary, F. \textit{On the notion of balance of a signed graph}. Michigan Math. J. 2. (1953-1954), 143-146 and addendum preceding p. 1.

\bibitem{kaufL} Kauffman, Louis H.; Lambropoulou, Sofia On the classification of rational tangles. Adv. in Appl. Math. 33 (2004), no. 2, 199-237.

\bibitem{kur}
Kuratowski, K., \textit{Sur le probleme des courbes gauches en topologie}, Fund. Math. 15 (1930), 271-283.
	
\bibitem{robertson seymour thomas}
	Robertson, N., Seymour, P.D, and Thomas, R.  \textit{Linkless embeddings of graphs in 3-Space}.  Bulletin (New Series) of the American Mathematical Society.  Vol 28 (1993), No. 1, 84-89.
	
\bibitem{sphere} Rukhovich, A., \textit{On intersection of two embedded spheres in 3-space}, Topology and its Applications
Volume 170, 15 June 2014, 96-103.
	
\bibitem{sachs} 
	Sachs, H. \textit{On spatial representations of finite graphs}. Finite and infinite sets, Vol. I, II (Eger, 1981), 649–662, Colloq. Math. Soc. János Bolyai, 37, North-Holland, Amsterdam, 1984. 
\bibitem{schu} Schubert, H., Knoten miot zwei Br\"ucken, Math. Zeitchriften, 65 (1956), 133-170.	

\bibitem{tutte}
Tutte, W. T. \textit{A homotopy theorem for matroids}. I, II. Trans. American Mathematical Society 88 (1958), 144-174.	
	
\bibitem{west}
West, D., Introduction to Graph Theory, Prentice Hall, Inc., Upper Saddle River, NJ, 1996.

\bibitem{wu}
	Wu, Ying-Qing.  \textit{On planarity of graphs in 3-manifolds}.  Commentarii Mathematici Helvetici.  67 (1992), no. 4, 635-647.  
	
	
	
	

\end{thebibliography}
\end{document}